\documentclass[11pt]{amsart}
\input{diagrams}

\usepackage{bibentry}
\usepackage{amsmath}
\usepackage{amsthm}
\usepackage{amssymb}
\usepackage{amsfonts}
\usepackage{amsxtra}
\usepackage{amscd}
\usepackage{epsfig}
\usepackage{verbatim}
\usepackage{latexsym,amstext,epsfig}

\usepackage[all, knot]{xy}
\xyoption{arc}

\newsymbol\pp 1275
\newcommand{\Hom}{\operatorname{Hom}}
\newcommand{\End}{\operatorname{End}}
\newcommand{\Ext}{\operatorname{Ext}}
\newcommand{\ext}{\operatorname{ext}}
\newcommand{\Rep}{\operatorname{Rep}}

\newcommand{\GL}{\operatorname{GL}}

\newcommand{\ZZ}{\mathbb Z}
\newcommand{\CC}{\mathbb C}

\newcommand{\coker}{\operatorname{coker}}

\newtheorem{theorem}{Theorem}[section]
\newtheorem{proposition}[theorem]{Proposition}

\theoremstyle{definition}
\newtheorem{definition}[theorem]{Definition}

\title[]{On orbit closures for infinite type quivers}

\author{Calin Chindris}
\address{University of Minnesota, School of Mathematics, Minneapolis, MN, USA}
\email[Calin Chindris]{chindris@math.umn.edu}

\date{August 24, 2007}
\begin{document}
\bibliographystyle{plain}
\begin{abstract}
For the Kronecker quiver, Zwara \cite[Theorem 1]{Zw2} has found an
example of a representation whose orbit closure is neither
unibranch nor Cohen-Macaulay. In this note, we explain how to
extend this example to all infinite type quivers without oriented
cycles.
\end{abstract}
\maketitle

In this short note, we use quiver exceptional sequences to reduce,
in a "hom-controlled" manner, the list of all
representation-infinite quivers without oriented cycles to just
$\widetilde{\mathbb{A}}_2$. This observation combined with Zwara's
example from \cite[Theorem 1]{Zw2} yields:

\begin{theorem}\label{main-thm} Let $Q$ be a representation-infinite, connected quiver without
oriented cycles. Then, there exists a representation $W$ whose
orbit closure is neither unibranch nor Cohen-Macaulay.
\end{theorem}

In \cite{Zw1}, Zwara showed that the orbit closures of
representations of Dynkin quivers are always unibranch. Hence, we
deduce that the Dynkin quivers are precisely the ones with the
property that all orbit closures are unibranch.

In what follows, we first review some background material from
quiver theory and then prove Theorem \ref{main-thm}. Throughout
this note, we work over an algebraically closed field $k$ of
characteristic zero.

Let $Q=(Q_0,Q_1,t,h)$ be a finite quiver, where $Q_0$ is the set
of vertices, $Q_1$ is the set of arrows and $t,h:Q_1 \to Q_0$
assign to each arrow $a \in Q_1$ its tail \emph{ta} and head
\emph{ha}, respectively. A representation $V$ of $Q$ over $k$ is a
family of finite dimensional $k$-vector spaces $\lbrace V(x) \mid
x\in Q_0\rbrace$ together with a family $\{ V(a):V(ta)\rightarrow
V(ha) \mid a \in Q_1 \}$ of $k$-linear maps. If $V$ is a
representation of $Q$, we define its dimension vector $\underline
d_V$ by $\underline d_V(x)=\dim_{k} V(x)$ for every $x\in Q_0$.
Thus the dimension vectors of representations of $Q$ lie in
$\Gamma=\ZZ^{Q_0}$, the set of all integer-valued functions on
$Q_0$.

Given two representations $V$ and $W$ of $Q$, we define a morphism
$\varphi:V \rightarrow W$ to be a collection of linear maps
$\lbrace \varphi(x):V(x)\rightarrow W(x)\mid x \in Q_0 \rbrace$
such that for every arrow $a\in Q_1$, we have
$\varphi(ha)V(a)=W(a)\varphi(ta)$. We denote by $\Hom_Q(V,W)$ the
$k$-vector space of all morphisms from $V$ to $W$. Let $W$ and $V$
be two representations of $Q.$ We say that $V$ is a
subrepresentation of $W$ if $V(x)$ is a subspace of $W(x)$ for all
vertices $x \in Q_0$ and $V(a)$ is the restriction of $W(a)$ to
$V(ta)$ for all arrows $a \in Q_1$. In this way, we obtain the
abelian category $\Rep(Q)$ of all quiver representations of $Q$.

A representation $W$ is said to be a \emph{Schur representation}
if $\End_{Q}(W) \cong \CC$. The dimension vector of a Schur
representation is called a \emph{Schur root}.

\textbf{From now on, we assume that our quivers are without
oriented cycles}. For two quiver representations $V$ and $W$,
consider Ringel's canonical exact sequence \cite{R} :
\begin{equation}\label{can-exact-seq}
0 \rightarrow \Hom_Q(V,W) \rightarrow \bigoplus_{x\in
Q_0}\Hom_{k}(V(x),W(x)){\buildrel
d^V_W\over\longrightarrow}\bigoplus_{a\in
Q_1}\Hom_{k}(V(ta),W(ha)),%\rightarrow \Ext^1(V,W)\rightarrow 0
\end{equation}
where $d^V_W((\varphi(x))_{x\in
Q_0})=(\varphi(ha)V(a)-W(a)\varphi(ta))_{a\in Q_1}$ and
$\Ext^1_{Q}(V,W)=\coker(d^V_W).$

If $\alpha,\beta$ are two elements of $\Gamma$, we define the
Euler inner product by
\begin{equation}\label{Euler-prod}
\langle\alpha,\beta \rangle_{Q} = \sum_{x \in Q_0}
\alpha(x)\beta(x)-\sum_{a \in Q_1} \alpha(ta)\beta(ha).
\end{equation}
(When no confusion arises, we drop the subscript $Q$.)

It follows from (\ref{can-exact-seq}) and (\ref{Euler-prod}) that
$$\langle\underline d_V, \underline d_W\rangle=
\dim_k\Hom_{Q}(V,W)-\dim_k\Ext^1_{Q}(V,W).$$

A dimension vector $\beta$ is called a \emph{real Schur} root if
there exists a representation $W \in \Rep(Q,\beta)$ such that
$\End_Q(W) \simeq k$ and $\Ext_{Q}^1(W,W)=0$ (we call such a
representation \emph{exceptional}). Note that if $\beta$ is a real
Schur root then there exists a unique, up to isomorphism,
exceptional $\beta$-dimensional representation.

For $\alpha$ and $\beta$ two dimension vectors, consider the
generic $\ext$ and $\hom$:
$$\ext_Q(\alpha,\beta):=\min \{\dim_k \Ext^1_Q(V,W) \mid (V,W)\in \Rep(Q,\alpha)\times \Rep(Q,\beta)\},$$
and
$$
\hom_Q(\alpha,\beta):=\min \{\dim_k \Hom_Q(V,W) \mid (V,W)\in
\Rep(Q,\alpha)\times \Rep(Q,\beta)\}.
$$
Given two dimension vectors $\alpha$ and $\beta$, we write $\alpha
\perp \beta$ provided that
$\ext_Q(\alpha,\beta)=\hom_Q(\alpha,\beta)=0$.

\begin{definition} We say that $(\varepsilon_1,\dots,\varepsilon_r)$ is an \emph{exceptional
sequence} if
\begin{enumerate}
\renewcommand{\theenumi}{\arabic{enumi}}

\item the $\varepsilon_i$ are real Schur roots;

\item $\varepsilon_i \perp \varepsilon_j$ for all $1 \leq i<j \leq
l$.

\end{enumerate}
\end{definition}

Following \cite{DW2}, a sequence
$(\varepsilon_1,\dots,\varepsilon_r)$ is called a \emph{quiver
exceptional sequence} if it is exceptional and $\langle
\varepsilon_j,\varepsilon_i \rangle \leq 0$ for all $1 \leq i<j
\leq l$.

Now, let $\varepsilon=(\varepsilon_1,\dots,\varepsilon_r)$ be a
quiver exceptional sequence and let $E_i \in
\Rep(Q,\varepsilon_i)$ be exceptional representations. Construct a
new quiver $Q(\varepsilon)$ with vertex set $\{1,\dots,r\}$ and
$-\langle \varepsilon_j,\varepsilon_i \rangle$ arrows from $j$ to
$i$. Define $\mathcal{C}(\varepsilon)$ to be the smallest full
subcategory of $\Rep(Q)$ which contains $E_1, \dots, E_r$ and is
closed under extensions, kernels of epimorphisms, and cokernels of
monomorphisms.

For the remaining of this section, we assume that $r \leq N-1$,
where $N$ is the number of vertices of $Q$. We recall a very
useful result from \cite[Section 2.7]{DW2} in a form that is
convenient for us (see also \cite{CC6}):

\begin{proposition}\cite{DW2} The category $\mathcal{C}(\varepsilon)$ is naturally equivalent
to $\Rep(Q(\varepsilon))$ with $E_1, \dots, E_r$ being the simple
objects of $\mathcal{C}(\varepsilon)$. Furthermore, the inverse
functor from $\Rep(Q(\varepsilon))$ to $\mathcal{C}(\varepsilon)$
is a full exact embedding into $\Rep(Q)$.
\end{proposition}

Consider the linear transformation
$$
\begin{aligned}
I: \ZZ^{Q(\varepsilon)_0}& \to \ZZ^{Q_0}\\
\alpha=(\alpha(1), \dots, \alpha(r))& \to \sum_{i=1}^r \alpha(i)
\varepsilon_i,
\end{aligned}
$$

If $F:\Rep(Q(\varepsilon)) \to \Rep(Q)$ is the full exact
embedding form Proposition \ref{excep-prop} then $F$ is clearly
hom-controlled in the sense of Zwara \cite{Zw3}. Let
$$
F^{(\alpha)}:\Rep(Q(\varepsilon), \alpha) \to \Rep(Q, I(\alpha))
$$
be the induced morphism of varieties. Zwara \cite[Theorem 2]{Zw3}
(see also \cite[Proposition 9]{BZ1}) showed that hom-controlled
functors preserve the type of singularities of orbit closures. In
particular, we have:

\begin{proposition}\label{excep-prop} Keep the same notations as above. Let $V \in \Rep(Q(\varepsilon),
\alpha)$. Then, $\overline{\GL(\alpha)V}$ is normal (or unibranch
or Cohen-Macaulay) if and only if
$\overline{\GL(I(\alpha))F^{(\alpha)}(V)}$ has the same property.
\end{proposition}

Next, let us recall Zwara's example from \cite{Zw2}:

\begin{theorem}\label{Zwara-thm} Let $\theta(2)$ be the Kronecker
quiver
$$
\xy     (0,0)*{1}="a";
        (10,0)*{2}="b";
        {\ar@2{->} "a";"b" };
\endxy
$$
Label the arrows by $a$ and $b$. Consider the following
representation $V \in \Rep(\theta(2),(3,3))$ defined by $V(a)=
\left(
\begin{matrix}
0& 0 & 0\\
1& 0 & 0\\
0& 1 & 0
\end{matrix}
\right)$ and $V(b)= \left(
\begin{matrix}
1& 0 & 0\\
0& 0 & 0\\
0& 0 & 1
\end{matrix}
\right)$. Then, $\overline{\GL(\alpha)V}$ is neither unibranch nor
Cohen-Macaulay.
\end{theorem}

\begin{proof}[Proof of Theorem \ref{main-thm}] From \cite[Lemma 2.1,
pp. 253]{ASS}, we know that any finite, connected quiver $Q$ of
infinite representation type must contain a Euclidean quiver as a
subquiver. Therefore, it is enough to prove the theorem for
Euclidean quivers. Now, let $Q$ be a Euclidean quiver and denote
by $\delta_Q$ the isotropic Schur root of $Q$.

Choose $v$ to be a vertex such that $Q \setminus v$ is a Dynkin
quiver. Without loss of generality, let us assume that $v$ is a
source. In this case, we take $\varepsilon_1=\delta_Q-e_v$ and
$\varepsilon_2=e_v$. Then, $(\varepsilon_1,\varepsilon_2)$ is a
quiver exceptional sequence with $\langle
\varepsilon_2,\varepsilon_1 \rangle=-2$ and so the proof follows
from Proposition \ref{excep-prop} and Theorem \ref{Zwara-thm}.
\end{proof}

%\bibliography{biblio}
\end{document}